\documentclass[12pt]{article}
\usepackage{amssymb}
\usepackage{amsthm}
\usepackage{latexsym,array}
\usepackage{amsthm}
\usepackage{amsmath, mathrsfs}
\usepackage[english]{babel}
\usepackage{cite}
\usepackage{color}

\newtheorem{theorem}{Theorem}[section]
\newtheorem{definition}[theorem]{Definition}

\newtheorem{proposition}[theorem]{Proposition}
\newtheorem{corollary}[theorem]{Corollary}
\newtheorem{remark}[theorem]{Remark}

\newcommand{\D}{\mathbb{D}}
\newcommand{\e}{{\mathbf{e}}}
\newcommand{\edag}{{\mathbf{e}^\dagger}}

\newcommand{\bc}{\mathbb{B}\mathbb{C}}

\newcommand{\R}{{\mathbb R}}
\newcommand{\C}{{\mathbb{C}}}

\newcommand{\ki}{{\bf k}}
\newcommand{\li}{{\bf i}}
\newcommand{\ji}{{\bf j}}
\begin{document}
\title{A bicomplex $(\vartheta,\varphi)-$weighted fractional Borel-Pompeiu type formula}
\small{
\author {Jos\'e Oscar Gonz\'alez-Cervantes$^{(1)}$ and Juan Bory-Reyes$^{(2)\footnote{corresponding author}}$}
\vskip 1truecm
\date{\small $^{(1)}$ Departamento de Matem\'aticas, ESFM-Instituto Polit\'ecnico Nacional. 07338, Ciudad M\'exico, M\'exico\\ Email: jogc200678@gmail.com\\$^{(2)}$ {SEPI, ESIME-Zacatenco-Instituto Polit\'ecnico Nacional. 07338, Ciudad M\'exico, M\'exico}\\Email: juanboryreyes@yahoo.com
}
\maketitle
\begin{abstract}
\small{
The purpose of this paper is to establish a Borel-Pompeiu type formula induced from a fractional bicomplex $(\vartheta,\varphi)-$weighted Cauchy-Riemann operator, where the weights are two hyperbolic orthogonal bicomplex functions and the fractionality is understand in the Riemann-Liouville sense.}
\end{abstract} 
\noindent
\textbf{Keywords.} Bicomplex analysis; Riemann-Liouville derivative; Cauchy-Riemann operator; Borel-Pompeiu formula.\\
\textbf{MSC Classification Numbers:} 26A33, 30A05, 30E20, 30G35, 32A30.

\section{Introduction and Preliminaries} 
Fractional calculus, involving derivatives-integrals of arbitrary real or complex order, is the natural generalization of the classical calculus, which in the latter years became a well-suited tool by many researchers working in several branches of science and engineering.

The interest in the subject has been growing continuously during the last few decades due to its wide range of theoretical and practical features and benefits, see \cite{GM, KST, OS, O, P, SKM} and the references given there. For a brief history and exposition of the foundations of the fractional calculus we refer the reader to \cite{MR, Ro}.

In 1892, Corrado Segre \cite{S} presented the bicomplex numbers system, see also \cite{Sco, Spa}. A deeper discussion of the historical appearance of bicomplex numbers, from their origin as \textit{imaginary} solutions of irrational equations, to their introduction in the study context of the algebras of hypercomplex numbers can be found in \cite{CC}.

The book of G. B. Price \cite{Pr} contains an evaluation of the subject of bicomplex numbers including a theory of differentiability in bicomplex numbers context and had to be regarded as the foundational work in this theory. For a recent account of the theory we refer the reader to \cite{EMDA}.

Fractional bicomplex calculus in the Riemann–Liouville sense is a very recent topic of research. The initial study has been presented in \cite{CTOP}. 

A Borel-Pompieu formula induced by a complex $\psi$-weighted Cauchy-Riemann operator appear in \cite{ATV}.

The aim of this paper is to obtain a Borel-Pompeiu type formula induced from a fractional bicomplex weighted Cauchy-Riemann operator, where the weights are two hyperbolic orthogonal bicomplex functions and the fractionality is understand in the Riemann-Liouville sense. Our approach will be different from \cite{CTOP} which has been done by the assumption of the one dimensional fractional derivative in the Riemann-Liouville sense in each direction over the bicomplex structure. 

\subsection{The Riemann-Liouville fractional derivatives and integrals}
There are different definitions of fractional derivatives. One of the most popular (even though it has disadvantages for applications to real world problems) is the Riemann–Liouville derivative (see, e.g., \cite{MR}). We begin by presenting a brief exposition of the basic ideas on Riemann-Liouville fractional derivatives and integrals.

Let us assume $a, b, \alpha \in \mathbb R$ such that $\alpha> 0$ and $a < b$. Suppose $f \in L^1([a, b], \mathbb R)$ identically equal to zero outside of the interval $[a, b]$.

The left and on the right Riemann-Liouville integrals of order $\alpha$ of $f$ are defined, respectively, as follows: 
$$({\bf I}_{a^+}^{\alpha} f)(x) := \frac{1}{\Gamma(\alpha)} \int_a^x \frac{f(\tau)}{(x-\tau)^{1-\alpha}} d\tau, \ \textrm{with} \ x > a$$
and
$$({\bf I}_{b^-}^{\alpha} f)(x) := \frac{1}{\Gamma(\alpha)} \int_x^b \frac{f(\tau)}{(\tau-x)^{1-\alpha}} d\tau, \ \textrm{with} \ x < b.$$
We will denote by $f\in AC^1([a, b], \mathbb R)$ the class of functions $f$ which are continuously differentiable and absolutely continuous on the segment $[a, b]$. The fractional derivatives in the Riemann-Liouville sense, on the left and on the right, are, respectively, defined by 
\begin{align}\label{FracDer} 
(D _{a^+}^{\alpha} f)(x):= \frac{d}{dx} \left[ ({\bf I}_{a^+}^{1-\alpha} f)(x)\right]
\end{align}
and
\begin{align} \label{FracDer1}
(D _{b^-}^{\alpha} f)(x):= (-1)\frac{d}{dx}\left[({\bf I}_{b^-}^{1-\alpha}f)(x)\right]. 
\end{align}
It is worth noting that the derivatives in (\ref{FracDer}), (\ref{FracDer1}) exist for $f\in AC^1([a, b], \mathbb R)$. Fractional Riemann-Liouville integral and derivative are linear operators.

Fundamental theorem for Riemann-Liouville fractional calculus \cite{CG} shows that 
\begin{align}\label{FundTheorem}
(D_{a^+}^{\alpha} {\bf I}_{a^+}^{\alpha}f)(x)=f(x) \  \textrm{and} \ (D _{b^-}^{\alpha}  {\bf I}_{b^-}^{\alpha} f)(x) = f(x).
\end{align}

Let us mention an important property of the fractional Riemann-Liouville derivative, see \cite[pag. 1835]{VTRMB}.
\begin{proposition}
\begin{equation}\label{cte}
( D _{a^+}^{\alpha} 1)(x)=\frac{(x-a)^{-\alpha}}{\Gamma[1-\alpha]}, \ \forall   x\in [a, b].
\end{equation}
\end{proposition}
\begin{remark}
Property (\ref{cte}) reflects the violation of Leibniz rule as a characteristic property of derivatives of non-integer orders, see \cite{Tarasov}.
\end{remark}
\subsection{A complex $\psi$-weighted Cauchy-Riemann operator} 
For convenience of the reader we repeat the relevant material from \cite{ATV} without proofs, thus making our exposition self-contained.

Given  $z,w\in \mathbb C$ is meant a product in the form 
\begin{align*}
\langle z   ,  w  \rangle_{\mathbb C}:= \frac{1}{2}\left(  \overline z  w + \overline w z \right)=
\frac{1}{2}\left(  z   \overline w +  w  \overline z \right).
\end{align*}
 
Consider functions $\psi_0 , \psi_1  \in  C^1(\C,\C)$  such that  $\langle \psi_0   ,  \psi_1   \rangle_{\mathbb C} = 0 $ on $\C$.
Let us introduce the notation $\psi:=(\psi_0,\psi_1 )$, where $\psi_0 = p_0  + ip_1$ and $\psi_1 =  q_0+ iq_1$. We thus get that  
$\langle \psi_0 ,  \psi_1   \rangle_{\mathbb C} \equiv  0$ on $\C$ if and only if  $p_ 0   \psi_1  = - i q_1\psi_0$  on $\C$. 

From the above, the following differential operator (to be called the $\psi$-weighted Cauchy-Riemann operator) is defined:
\begin{align*}
\mathcal D_{\psi} :=  \psi_0 \frac{\partial}{\partial x }  +   \psi_1 \frac{\partial}{\partial y }.   
\end{align*}

\begin{definition}
A continuously differentiable complex-valued function $f$ in a regular domain (continuous differentiable or smooth on the boundary) $\Omega \subset \mathbb C$ is said to be $\psi$-weighted analytic in $\Omega$ if $f$ satisfies the equation $\mathcal D_{\psi} f = 0$ in $\Omega$.
\end{definition}
Holomorphic functions are special cases of null solutions of $\mathcal D_{\psi}$, when $\psi_0=1 $ and $ \psi_1 = {i}$.  
\begin{theorem} ($\psi$ -weighted Gauss theorem) Assume that $\Omega \subset \C$ is a regular domain and let $f$, $\psi_0 = p_0 + { i}p_1$, $\psi_1 = q_0 + {i}q_1$ in  $C^1(\Omega, \C) \cap  C(\overline \Omega, \C)$. Then
\begin{align*}  
\int_\Omega\left( 
{\mathcal D}_{\psi} f(w) + ( \frac{\partial p_0}{\partial x }   +  \frac{\partial q_0}{\partial y }  )f + ( \frac{\partial p_1}{\partial x }  +  \frac{\partial q_1}{\partial y }  ) {  i} f
\right)
dxdy =
\int_{\partial \Omega}
 f d\sigma_{\psi}(w), 
\end{align*}
where $w=x+iy$ and the differential form $d\sigma_{\psi}(w):= \psi_0(w)dy - \psi_1(w)dx.$
\end{theorem}

\begin{theorem} 
($\psi$ -weighted Cauchy-Pompeiu integral formula) \label{weighted Cauchy-Pompeiu}
Let $\Omega \subset \C$ be a regular domain, $f \in  C^1(\Omega,\C)\cap C(\Omega,\C)$ and $\psi_0, \psi_1\in \C$. For any $z\in \Omega$, we have the Cauchy-Pompeiu integral
formula
\begin{align*}  f(z) c_{\psi} =
\int_{\partial \Omega}  
f( w )E_{\psi}( w , z)d\sigma_{\psi}( w ) - \int_\Omega 
E_{\psi}( w ,z) {\mathcal D}_{\psi} f( w )dx d y,
\end{align*}
 where $ w  = x + {  i} y$,  
\begin{align*}
c_{\psi} = \sin^2(\alpha) \int_0^{2\pi}
\frac{\psi_0
     \cos \theta + \psi_1 \sin \theta}{
\left( \dfrac{\cos \theta}
{r_0^2} - \dfrac{
\sin \theta \cos \alpha}{r_0r_1}
\right)
 \psi_0 +
\left(
\dfrac{\sin \theta}
{r^2_1} - 
 \dfrac{ \cos \theta \cos \alpha}{
r_0r_1}
\right)
 \psi_1}
d\theta
\end{align*}
and $E_\psi( \cdot, \cdot )$ is the $\psi$-Cauchy-type kernel, see \cite[Section 3]{ATV}.  
\end{theorem} 
For completeness, we will discuss further information about the weighted complex analytic function theory.

Let $\psi_0,\psi_1\in \mathbb C\setminus\{0\}$ orthogonal complex numbers. It is easily seen that there exist $a,b,c,d\in \mathbb R$ such that  $a \psi_0 + b\psi_1 =1$ and $c \psi_0 + d\psi_1 =i$.

\begin{enumerate}
\item We introduce the following real-linear transformation 
$$T(z')=T(x',y') = \left(\begin{array}{ cc} a & b \\ c & d \end{array} \right)  \left( \begin{array}{c}x' \\ y' \end{array} \right) = \left( \begin{array}{c}x \\ y \end{array} \right)=T(z),$$
where $z'=x'+iy'$ and  $z= x+ i y$. For $f \in C^1(\Omega, \mathbb C)$, where $\Omega\subset\mathbb C$ is a domain of complex-variable $z$, it follows that
\begin{align}\label{RelaComplex1}
 \mathcal D_{\psi, z'}(f\circ T)   = &  \psi_0 \frac{\partial f\circ T}{\partial x' }   +   \psi_1 \frac{\partial f\circ T}{\partial y' }   \nonumber \\
                           = &  \psi_0 a \frac{\partial f}{\partial x } \circ T   + \psi_0 c \frac{\partial f}{\partial y } \circ T 
													 +  \psi_1 b  \frac{\partial f}{\partial x } \circ T   + \psi_1 d \frac{\partial f}{\partial y } \circ T  \nonumber \\
													= &  (a\psi_0    +  b\psi_1  ) \frac{\partial f}{\partial x } \circ T 
													+  (c \psi_0    +  d\psi_1  )  \frac{\partial f}{\partial y } \circ T   \nonumber \\
																				= &    \frac{\partial f}{\partial x } \circ T 
													+ i \frac{\partial f}{\partial y } \circ T = \left( \frac{d  }{ d \bar z} f\right) \circ T, 
\end{align}
where $\mathcal D_{\psi, z'}$ represents the weighted Cauchy-Riemann operator $\mathcal D_{\psi}$ with respect to the real components of $z'$. 

In other words,  
\begin{align}\label{RelaComplex2}
 \mathcal D_{\psi, z'}(g) = \left(\frac{d  }{ d \bar z}  g\circ T^{-1} \right)\circ T. 
\end{align}
for every function $g$ of class $C^1$ on a domain of the complex-variable $z'$.
\item \textit{A  Gleason type problem for weighted analytic functions}. Let us denote by $\mathbb D(0,1)$ the unit disk in the complex plane $\C$.  If $f\in Ker \mathcal D_{\psi} \cap C^{1}( \mathbb D(0,1),\mathbb C)$, Theorem \ref{weighted Cauchy-Pompeiu} yields   
$f\in  C^{\infty}(\mathbb D(0,1),\mathbb C)$ and the fundamental theorem of calculus applied to each real component of $f$ leads to  
\begin{align*}
f(z)-f(0) &= \int_{0}^1 \frac{d f(tz)}{d t} dt  \\
 &=   \int_{0}^1 \left( \frac{\partial  f }{\partial x}(tz)  x  +  \frac{\partial f}{\partial y}(tz) y \right) dt. 
\end{align*}
As  
\begin{align*}
\mathcal D_{\psi}f(tz)  =  \psi_0 \frac{\partial f}{\partial x } (tz)  +   \psi_1 \frac{\partial f}{\partial y } (tz) =0   
\end{align*}
we have
\begin{align*}
  \frac{\partial f}{\partial y } (tz) =- \frac{ \psi_0}{\psi_1} \frac{\partial f}{\partial x } (tz).       
\end{align*}
Therefore,  
\begin{align*}
f(z)-f(0) &=    \int_{0}^1 \left( \frac{\partial  f }{\partial x}(tz)  x    - \frac{ \psi_0}{\psi_1}  y \frac{\partial f}{\partial x } (tz)  \right) dt \\
 &  =    \int_{0}^1 \left( x    - \frac{ \psi_0}{\psi_1}  y  \right)  \frac{\partial  f }{\partial x}(tz)  dt 
 = \left( x    - \frac{ \psi_0}{\psi_1}  y  \right)  h_1(z),\\ 
\end{align*}
where 
$$h_1(z) = \int_{0}^1  \frac{\partial  f }{\partial x}(tz)  dt.$$
The Leibniz rule allows us to see that $h_1\in  Ker \mathcal D_{\psi} \cap C^{\infty}( \mathbb D(0,1),\mathbb C)$. Repeating previous computation for $h_1$ enables us to write  
\begin{align*}
f(z)= a_0 +  \left( x - \frac{ \psi_0}{\psi_1}  y  \right)   a_1 +   \left( x    - \frac{ \psi_0}{\psi_1}  y  \right)^2  h_2(z),\\
\end{align*}
where
$$h_2(z) =\int_0^1 \frac{\partial    }{\partial x_2} \int_{0}^1  \frac{\partial  f }{\partial x_1}(t_1 t_2z)  dt_1 dt_2= \int_0^1 \int_{0}^1 \frac{\partial    }{\partial x_2}  \frac{\partial  f }{\partial x_1}(t_1 t_2z)  dt _1dt_2.$$
We continue successively in this fashion to prove the existence of a sequence of complex numbers $(a_n) $ such that  
 \begin{align*}
f(z)= & \sum_{n=0}^{\infty}\left( x    - \frac{ \psi_0}{\psi_1}  y  \right)^n   a_n  
=   \sum_{n=0}^{\infty}\left( \frac{z +\bar z}{2}    - \frac{ \psi_0}{\psi_1}  \frac{z -\bar z}{2}  \right)^n   a_n  \\
=  & \sum_{n=0}^{\infty}\left( (1 - \frac{ \psi_0}{\psi_1} )  z +  (1 + \frac{ \psi_0}{\psi_1} )   \bar z  \right)^n  \frac{ a_n}{2^n},    
\end{align*}
where the convergence is uniform on compacts subsets of $\mathbb D(0,1)$. 

In summary,  every $f\in  Ker \mathcal D_{\psi} \cap C^{1}( \mathbb D(0,1),\mathbb C)$ satisfies certain version of Taylor's Theorem.
\item Let $a\in \mathbb C$ and consider $|w-a|=r$ and $ |z-a|<r$. The usual Cauchy kernel can be approximated as follows  
 \begin{align*}
K(w-z) = & \frac{1}{2\pi i} \frac{1}{w-z} = \frac{1}{2\pi i} \frac{1}{w-z} = \frac{1}{2\pi i (w-a)} \frac{1}{1-\frac{z-a}{w-a}} \\
\\ 
 = &  \frac{1}{2\pi i } \sum_{n=0}^{\infty} (z-a)^n(w-a)^{-(n+1)} =   \frac{1}{2\pi i } \sum_{n=0}^{\infty} u_n(z)v_n(w), 
\end{align*}  
where $u_n(z)= (z-a)^n$,  $v_n(w)= (w-a)^{-(n+1)}$ and the series is normal convergent.

The identities \eqref{RelaComplex1} and \eqref{RelaComplex2} help us to extends crucial ingredients of the usual holomorphic function theory to that induced by $ \mathcal D_{\psi, z'}.$  For instance, the weighted Cauchy kernel can be expressed in the form  
\begin{align*}
E_{\psi}( w' - z' ) = K (T(w')-T(z')) = K\circ T (w'- z')  =  \frac{1}{2\pi i } \sum_{n=0}^{\infty}U_n(z')V_n(w'), 
\end{align*}  
where $U_n = u_n\circ T $, \ $V_n =v_n\circ T $ and the series is normal convergent. 
\end{enumerate}
\subsection{Basic definitions and tools of bicomplex function theory}

The set $\bc$ of bicomplex numbers is defined as: 
$$\bc := \{Z=z_{1}+ z_{2}\ji;\, z_{1}, z_{2}\in {\C}\},$$
where $\C$ is the set of complex numbers with the imaginary unit $\li$. Here and subsequently, $\li$ and $\ji \neq \li$ stand for commuting imaginary units such that 
\begin{equation}\label{unidadesbico}
\li^{2}=\ji^{2}=-1; \li\ji=\ji\li.
\end{equation}
Let us denote by $\C(\li) := \{x + \li y : x, y \in {\R}\}$  and $\C(\ji) := \{x + \ji y : x, y \in {\R}\}$. Clearly $\C(\li)$ and $\C(\ji)$ are isomorphic to $\C$, but they co-exist inside $\bc$ as two different sets. Throughout the paper $\C$ is understood as $\C(\li)$. 

Let $\{\e, \edag \}$ be the idempotent base for $\bc$ over $\C(\li)$ and $\C(\ji)$, given by 
\begin{displaymath}
\e = \frac{1}{2}(1 + \li\ji), \hspace{5mm} \edag = \frac{1}{2}(1 - \li\ji).
\end{displaymath}
These idempotent elements satisfy that  
$\e \edag = 0$, 
$\e^2 = \e$, $(\edag
)^2 = \edag$;
$ \e + \edag
= 1$, $\e - \edag= {\bf ij}$.

Each bicomplex number may be written as: 
$$Z=z_1\e+z_2\edag,$$
where $z_1,z_2\in\C$. The latter is called the idempotent representation of $Z$, which has numerous advantages for calculations. For example, the product and the sum of bicomplex numbers are calculated component-wise. Specifically, if $Z =z_1 \e + z_2 \edag $ and $W =w_1 \e + w_2 \edag $ are two bicomplex numbers, then
$$Z +W=(z_1+w_1) \e + (z_2+w_2) \edag , \ Z W=(z_1w_1) \e + (z_2w_2) \edag .$$
The set of hyperbolic numbers
\[ \D := \{\lambda_1 + {\bf k} \lambda_2 \ \mid \ \lambda_1, \lambda_2 \in \mathbb R \},\]  
inherits all the algebraic definitions, operations and properties from $\bc$. 

For two elements $X,Y\in \bc$ we will write $X\preceq Y$  if and only if $Y-X\in \D^{+}$, where  
\begin{displaymath}
\D^{+} := \{ \lambda_1\e + \lambda_2\edag \ |\ \lambda_1 \in \R^{+}\ \land\ \lambda_2 \in \R^{+}\}.
\end{displaymath}
This relation define a partial order over $\bc$, see for instance \cite{EMDA}.

Given $Z=z_1\e +z_2\edag \in \bc$, we define the hyperbolic modulus of $Z$ as 
$$|Z|_{\ki}:=|z_1|\e + |z_2|\edag.$$

The hyperbolic-valued modulus induces the topology in $\bc$-module of the bicomplex balls 
$$B(W, r) := \{Z \in \bc \ \mid \   \  |Z-W|_{\bf k}   \prec r \}$$
with hyperbolic radius $r$ centered at $W\in \bc$. This is equivalent to the one induced by Euclidean norm under identification $\bc\cong \R^{4}$.

The set $\bc$ admits some other forms of writing, which help us to visualize its structure. 
\begin{align*} 
\bc  = &  \{w_1 + { \bf j} w_2 \ \mid \ z_1, z_2 \in \C({\bf i}) \} \\ 
   = &   \left\{    \mathfrak z_1  + {\bf i} \mathfrak z_2 \  \mid \ \mathfrak z_1 ,    \mathfrak z_2 \in \D \right\} \\
= &  \{z_1 \e +  z_2\edag  \ \mid \ z_1, z_2 \in \C({\bf i}) \}.
\end{align*}
Given $Z  = z_1 \e +  z_2\edag \in \bc$ introduce the conjugation $Z^* =  \overline z_1 \e +  \overline  z_2\edag$, where $\overline z_1, \overline z_2$ are usual complex conjugates to $z_1, z_2 \in \C({\bf i})$.

Note that 
$$ZZ^*=Z^*Z=  |Z|^2_{\bf k}.$$

Let $Z = z_1\e +z_2 \edag $ and  $W = w_1\e +w_2 \edag \in \bc$, we have 
\begin{align*} \langle Z, W\rangle_{{\bf k}} & := \frac{1}{2}\left({Z}^* W +  {W}^* Z   \right) = \frac{1}{2}\left( Z {W}^*  + W {Z}^*  \right)  \\
& =\frac{1}{2} ( \overline{z}_1 w_1 + \overline{w}_1z_1 ) \e  + 
\frac{1}{2}( \overline{z}_2 w_2 + \overline{w}_2z_2 ) \edag \\
& = \langle z_1,  w_1\rangle_{\C({\bf i })}  \e  + 
\langle z_2,  w_2 \rangle_{\C({\bf i })}  \edag. \\
\end{align*} 

\begin{definition}
Let $Z \in \Omega\subset\mathbb{BC}$ and $F: \Omega\rightarrow\mathbb{BC}$. The following limit, if exists, is called the derivative of $F$ at the point $Z$:
$$F'(Z):=\lim_{\Omega \ni W \to Z} \frac{F(W)-F(Z)}{W-Z},$$
such that $W-Z$ is an invertible bicomplex number.		
\end{definition}

\begin{definition}
A function $F:\Omega\subset\mathbb{BC}\rightarrow\mathbb{BC}$ is said to be $\mathbb{BC}$-holomorphic in $\Omega$ if for every $Z\in \Omega$ there exists the derivative $F'(Z)$.
\end{definition}

By \cite[Theorem 7.6.4]{EMDA}, a bicomplex valued function $F = f_1\e + f_2\edag $ defined on a product-type domain, i.e. $\Omega=\Omega_1 \e + \Omega_2 \edag \subset \bc$, where $\Omega_1 , \Omega_2  \subset \C({\bf i})$ are domains, is $\bc$-holomorphic if and only if  
$$F(Z) = f_1(z_1)\e + f_2(z_2)\edag$$
at every $Z=z_1{\e} +z_2\edag \in \Omega$, where $z_l\in \Omega_l$ and  $f_l\in C^1(\Omega_l, \C({\bf i})) \cap Ker \dfrac{d}{d \overline z_l} $ for $l=1,2$.

We follow the notation of \cite[Subsection 11.2]{EMDA} in assuming $\Omega$ to be a domain in $\bc$, and consider a two-dimensional, simply connected, piecewise smooth surface $\Gamma \in \Omega$ with boundary $\gamma= \partial \Gamma \subset \Omega$ which has the following properties: $\Gamma$ has a parametrization $R = R(u, v)$ such that $R = R_1\e+R_2\edag$, being $R_1$ and $R_2$ parametrization, respectively, of simply connected domains $\Gamma_1$ and $\Gamma_2$ in $\C({\bf i})$. Moreover, $\gamma$ is parametrized by $r = r(t)$, which is the restriction of $R$ onto $\partial \Gamma$ and is such that $r = r_1\e + r_2\edag$ where $r_1$ and $r_2$ are a parametrization, respectively, of $\gamma_1 := \partial \Gamma_1$ and of $\gamma_2 := \partial \Gamma_2$. The curves $\gamma_1$ and $\gamma_2$ are piecewise smooth, closed, Jordan curves in $\C({\bf i})$.

Let $F(Z) = f_1(z_1)\e+f_2(z_2)\edag$ a continuous bicomplex function on $\gamma$, where $f_l \in C^{1}(\Gamma_l, \C), l=1,2$, then a bicomplex integration is defined to be  
\begin{align*}
\int_{\gamma} F(Z) dZ := \e \int_{\gamma_1} f_1(z_1) d z_1 + \edag \int_{\gamma_2} f_2(z_2) d z_2, 
\end{align*} 
where $Z= z_1\e + z_2\edag $ and  $z_1,z_2\in \mathbb C({\bf i})$.

Similarly, consider
\begin{align*}
\int_{\Gamma} F(Z) dZ\wedge dZ^* :=
\left( \int_{\Gamma_1}f_1(z_1) d z_1 \wedge d\bar z_1
\right)\e   +\left( \int_{\Gamma_2} f_2(z_2) d z_2 \wedge d\overline z_2
\right)\edag . 
\end{align*} 
 
\begin{theorem}\label{BBPF}\cite{EMDA}(Bicomplex Borel–Pompeiu formula). Let $F \in C^1(\Omega, \bc)$ such that $F(Z) = f_1(z_1)\e+ f_2(z_2)\edag, \  Z = z_1\e+z_2\edag$, and let $\gamma$ and $\Gamma$ satisfy the above assumptions. Then for any $W \in \Gamma \setminus \gamma$, we have
\begin{align*}
F(W) = \frac{1}{2\pi{\bf i}} \int_{\gamma} \frac{F(Z)}{Z - W} dZ + \frac{1}{2\pi{\bf i}}\int_{\Gamma} \frac{\dfrac{\partial F}{\partial Z^*}}{Z - W}dZ \wedge dZ^*,
\end{align*}
where $Z = z_1\e + z_2\edag $ and  $z_1,z_2\in \mathbb C({\bf i})$ and
$\dfrac{\partial}{\partial Z^*} = \e \dfrac{\partial}{\partial \overline z_1} + \edag \dfrac{\partial}{\partial \overline z_2}.$
\end{theorem}
\begin{remark}
For a deeper discussion of bicomplex integration we refer the reader to \cite [Section 2]{BPS} and the references given there. 
\end{remark}
\section{A bicomplex $(\vartheta,\varphi)-$weighted holomorphic\\ functions}
Consider the functions $\vartheta (Z)= \vartheta_1(z_1) \e +  \vartheta_2(z_2)\edag $ and $ \varphi(Z) = \varphi_1(z_1)\e +\varphi_2(z_2)\edag$ 
of variable $Z=z_1\e +z_2\edag\in \bc$ such that $\vartheta_1, \vartheta_2, \varphi_1, \varphi_2\in  C^1(\C({\bf i}),\C({\bf i}))$. Then $\langle \vartheta  ,  \varphi \rangle_{\bf k} \equiv  0 $  on  $\C({\bf i})$ if and only if $\langle \vartheta_{\ell} ,   \varphi_{\ell}  \rangle_{\C({\bf i })} \equiv 0$ on $\bc$ for $\ell=1,2$. 

\begin{proposition}
	Denote 
$	\vartheta_1 =p_{_{1,1}}   + {\bf i} p_{_{1,2}} $, \  $\vartheta_2    =  p_{_{2,1}}  + {\bf i}p_{_{2,2}}$, \ 
 $ \varphi_1 = q_{_{1,1}}  + {\bf i} q_{_{1,2}}$ and $\varphi_2  =  q_{_{2,1}} + {\bf i} q_{_{2,2}}  $. 
Then    
$\langle \vartheta  ,  \varphi  \rangle_{\bf k} = 0$ on  $\bc$ 
if and only if 
\begin{align*}
 (p_{_{1,2}} \e   + p_{_{2,2}} \edag)    \varphi 
= - {\bf i} \left(  (q_{_{1,1}} \e     +   q_{_{2,1}} \edag\right)  \vartheta  ,   
\textrm{on} \ \bc.
\end{align*}
\end{proposition}
\begin{proof}
$\langle \vartheta  ,  \varphi  \rangle_{\bf k} = 0$ on  $\bc$ if and only if 
\begin{align*}\langle \vartheta_l  ,  \varphi_l  \rangle_{\C} = 0 \ \
\textrm{on}\  \C \ \textrm{for}\  l=1,2.
 \end{align*}
if and only if  
\begin{align*}
 p_{_{l,2}}   \varphi_l  = -{\bf  i}  q_{_{l,1}}  \vartheta _l  \  
\textrm{on}\ \C \ \textrm{for}\  l=1,2.
 \end{align*}
The following three statements are equivalent:
\begin{align*}
 p_{_{1,2}} \  \varphi_1 \e 
+ p_{_{2,2}} \  \varphi_2 \edag
= - {\bf i} \left(  q_{_{1,1}} \  \vartheta _1 \e
+  q_{_{2,1}} \  \vartheta _2  \edag \right)    
\textrm{ on } \ \bc,
 \end{align*}
\begin{align*}
 (p_{_{1,2}} \e   + p_{_{2,2}} \edag)   \   (    \varphi_1 \e 
+ \  \varphi_2 \edag ) 
= - {\bf i} \left(  (q_{_{1,1}} \e     +   q_{_{2,1}} \edag) \    ( \vartheta _1 \e
+   \vartheta _2  \edag )  \right)   
\textrm{ on } \ \bc,
 \end{align*}
\begin{align*}
 (p_{_{1,2}} \e   + p_{_{2,2}} \edag)    \varphi 
= - {\bf i} \left(  (q_{_{1,1}} \e     +   q_{_{2,1}} \edag\right)  \vartheta \    
\textrm{on} \ \bc,
 \end{align*}
which is clear from \cite[Proposition 2.1]{ATV}.
\end{proof}

\begin{definition}
Let $Z = z_1\e +z_2 \edag =  (x_1+ y_1{\bf  i })\e +(x_2+ y_2 {\bf i})\edag$ be a bicomplex variable and define   
\begin{align*}
 \frac{\partial}{\partial Z_{\vartheta \varphi } } := &     (\vartheta_1 \frac{\partial  }{\partial {x_1} } +  \varphi_1  \frac{\partial  }{\partial {y_1} }) \e  +  ( \vartheta_2 \frac{\partial  }{\partial {x_2} } +  \varphi_2 \frac{\partial  }{\partial {y_2} }) \edag.
\end{align*}
Explicitly, we can write   
\begin{align*} 
\frac{\partial  }{\partial Z_{\vartheta  \varphi } }  =  & \left(  ( p_{_{1,1}}  \frac{\partial  }{\partial {x_1} } 
 + q_{_{1,1}}   \frac{\partial  }{\partial {y_1} } )
+ {\bf i} (  p_{_{1,2}}    \frac{\partial  }{\partial {x_1} } +
     q_{_{1,2}}  \frac{\partial  }{\partial {y_1} }) \right)  \e   \\ 
		 &  + 
	\left(    ( p_{_{2,1}}   \frac{\partial  }{\partial {x_2} }  + q_{_{2,1}}  \frac{\partial  }{\partial {y_2} })  + {\bf i}  (p_{_{2,2}}   \frac{\partial  }{\partial {x_2} }  + q_{_{2,2}}  \frac{\partial  }{\partial {y_2} })  \right) \edag.
\end{align*}
We call $\dfrac{\partial  }{\partial Z_{\vartheta \varphi } }$ the bicomplex $(\vartheta,\varphi)-$weighted Cauchy-Riemann operator.
\end{definition} 
\begin{remark}
Let $\Omega \subset\bc$ a domain and let $F\in C^1(\Omega, \bc)$. There exist $f_1,f_2,f_3,f_4\in C^1(\Omega, \mathbb R)$ such that  $F= ( f_1\vartheta_1 + f_2\varphi_1) \e + ( f_3 \vartheta_2+ f_4\varphi_2)\edag.$ Hence $\dfrac{\partial  F}{\partial Z_{\vartheta \varphi } }=0$  on $\Omega$ if and only if 
$$ \vartheta_1^2 \frac{\partial f_1}{\partial {x_1}}  +  \varphi_1^2 \frac{\partial  f_2}{\partial {y_1}} + \vartheta_1 \varphi_1 ( \frac{\partial f_1}{\partial {y_1}}   +  \frac{\partial f_2 }{\partial {x_1}})=0,$$
$$ \vartheta_2^2 \frac{\partial f_3}{\partial {x_2}}  +  \varphi_2^2 \frac{\partial  f_4}{\partial {y_2}} + \vartheta_2 \varphi_2 ( \frac{\partial f_3}{\partial {y_2}}   +  \frac{\partial f_4 }{\partial {x_2}})=0$$
on $\Omega$.

On the other hand, following notation used in \cite{ATV} we have 
\begin{align*} 
\dfrac{\partial}{\partial Z_{\vartheta \varphi } } = &  D_{z_1, \vartheta_1 , \varphi_1 }  \e  +  D_{z_2, \vartheta_2 , \varphi_2 } \edag \\
 = &  \frac{1}{2}( D_{z_1, \vartheta_1 , \varphi_1 }   +  D_{z_2, \vartheta_2 , \varphi_2 }  ) + \frac{1}{2}( D_{z_1, \vartheta_1 , \varphi_1 } -  D_{z_2, \vartheta_2 , \varphi_2 }  ){\bf k}.
\end{align*}
Let $G=g_1+ g_2 {\bf k}$, where $g_1, g_2 \in C^1(\Omega\subset\bc, \C({\bf i}))$, then $\dfrac{\partial  G }{\partial Z_{\vartheta \varphi } }=0$ if and only if 
$$ D_{z_1, \vartheta_1 , \varphi_1 } g_1  +  D_{z_2, \vartheta_2 , \varphi_2 } g_1  +   D_{z_1, \vartheta_1 , \varphi_1 } g_2  -  D_{z_2, \vartheta_2 , \varphi_2 } g_2=0,$$
$$D_{z_1, \vartheta_1 , \varphi_1 } g_2  +  D_{z_2, \vartheta_2 , \varphi_2 } g_2 +   D_{z_1, \vartheta_1 , \varphi_1 }g_1   -  D_{z_2, \vartheta_2 , \varphi_2 } g_1=0.$$
\end{remark} 

\begin{theorem} \label{weighted_Gauss_Bicomplex}($(\vartheta,\varphi)-$ weighted Gauss theorem for bicomplex functions) 
Let $\Omega= \Omega_1\e+\Omega_2 \edag \subset \bc$ such that  $\Omega_1,\Omega_2\subset \mathbb C({\bf i})$ are domains 
and set $F =f_1 \e + f_2 \edag$ a bicomplex function with $f_l\in C^1(\Omega_l, \C) \cap  C(\overline \Omega_l, \C) $ for $l=1,2$.
Suppose a two-dimensional, simply connected, piecewise smooth surface $\Gamma \subset \Omega$ with smooth boundary $\gamma= \partial \Gamma \subset \Omega$ with properties according to Theorem \ref{BBPF}.  Then 
\begin{equation*}   
 \int_{\Gamma} \left(  \frac{\partial F }{\partial Z_{\vartheta  \varphi } } 
 +  A_{\vartheta  \varphi }  F  + B_{\vartheta  \varphi }  {\bf i} F  \right)
dZ \wedge dZ^{*}  = \int_{\gamma} F(Z) d\sigma_{\vartheta  \varphi }(Z) ,
\end{equation*}
where   $d\sigma_{\vartheta  \varphi }(Z) =  d\sigma_{\vartheta }(z_1) \e + d\sigma_{\varphi }(z_1) \edag$ and
\begin{align*}
A_{\vartheta  \varphi } =&   (\frac{\partial  p_{_{ 1,1}}  }{\partial  x_1} +\frac{\partial  q_{_{ 1,1}} }{\partial  y_1} ) \e  + (\frac{\partial  p_{_{ 2,1}}  }{\partial  x_2}+\frac{ \partial  q_{_{ 2,1}}}{\partial  y_2} ) \edag  
\\
B_{\vartheta  \varphi } = &  (
\frac{\partial  p_{_{ 1,2}} }{\partial  x_1} + \frac{\partial  q_{_{ 1,2}} }{\partial  y_1}   ) \e + 
 ( \frac{\partial  p_{_{ 2,2}}  }{\partial  x_2}  + \frac{\partial  q_{_{ 2,2}}  }{\partial  y_2}
) \edag  \end{align*}
\end{theorem}
\begin{proof}
A direct computation shows that  
\begin{align*}  &  \int_{\Gamma_1}
\left( 
{\mathcal D}_{\vartheta} f_1(z_1) +(\frac{\partial  p_{_{ 1,1}}  }{\partial  x_1} +\frac{\partial  q_{_{ 1,1}} }{\partial  y_1} ) f_1 +
(
\frac{\partial  p_{_{ 1,2}} }{\partial  x_1} + \frac{\partial  q_{_{ 1,2}} }{\partial  y_1}   ) {\bf i} f_1
\right)
dx_1dy_1  \e  \\
& + \int_{\Gamma_2}
\left( 
{\mathcal D}_{\varphi} f_2(z_2) + 
(\frac{\partial  p_{_{ 2,1}}  }{\partial  x_2}+\frac{ \partial  q_{_{ 2,1}}}{\partial  y_2} ) f_2 
 + ( \frac{\partial  p_{_{ 2,2}}  }{\partial  x_2}  + \frac{\partial  q_{_{ 2,2}}  }{\partial  y_2}
) {\bf i} f_2
\right)
dx_2dy_2  \edag\\
 = & 
\int_{\gamma_1}
f_1 d\sigma(z_1) \e +  
\int_{\gamma_2}f_2 d\sigma(z_2) \edag,  
\end{align*}
which proves the theorem.
\end{proof}
 
\begin{theorem}\label{weighted_Borel-Pompieu_Bicomplex} 
(Bicomplex $(\vartheta,\varphi)-$Borel-Pompieu formula) Suppose the hypothesis and notation of the previous theorem. Given $W\in \Gamma$ we get  
\begin{align*}  & F(W)c_{(\vartheta,\varphi)}  =
\int_{\gamma  }  
F (Z )E_{(\vartheta, \varphi)} (Z , W )d\sigma_{(\vartheta, \varphi)}(Z ) 
 - \int_{\Gamma} 
E_{(\vartheta , \varphi)}(Z ,W ) \frac{\partial F }{\partial Z_{(\vartheta,\varphi)}}  dZ\wedge dZ^{*},
\end{align*}
where $c_{(\vartheta,\varphi)} = c_{\vartheta_1 \varphi_1 }  \e +  c_{\vartheta_2 \varphi_2}\edag,$ 
$$E_{(\vartheta  , \varphi)} (Z,W) =  E_{\vartheta_1\varphi_1 } (z_1 , w_1 ) \e  +
 E_{\vartheta_2\varphi_2 } (z_2 , w_2)\edag$$ 
and 
$$d\sigma_{(\vartheta , \varphi)}(Z) =  d\sigma_{\vartheta_1\varphi_1 }(z_1 ) \e  + d\sigma_{\vartheta_2 \varphi_2 }(z_2 ) \edag.$$
\end{theorem}
\begin{proof}
Let $\Omega_l \subset \C$ be regular domains and $f_l \in  C^1(\Omega_l,C)\cap C(\Omega_l,\C)$, $l=1,2$. Given $w_l\in \Gamma_l$, the complex $(\vartheta_l , \varphi_l )$-weighted Cauchy-Pompeiu integral formula gives us that
\begin{align*}  & f_1 (w_1 ) c_{\vartheta_1 \varphi_1 }  \e +   f_2 (w_2 ) c_{\vartheta_2 \varphi_2 }  \edag  \\
 = & 
\int_{\gamma_1 }  
f_1 (z_1 )E_{\vartheta_1\varphi_1 } (z_1 , w_1 )d\sigma_{\vartheta_1 \varphi_1 }(z_1 ) \e  +
\int_{\gamma_2 }  
f_2 (z_2 )E_{\vartheta_2 \varphi_2 } (z_2 , w_2 )d\sigma_{\vartheta_2 \varphi_2 }(z_2 ) \edag \\
 & \  \ 
 - \int_{\Gamma_1 } 
E_{\vartheta_1 \varphi_1 }(z_1 ,w_1 ) D_{z_1, \vartheta_1 , \varphi_1 } f_1 (z_1 )dx_1  dy_1 \e  \\
  & \ \ + 
\int_{\Gamma_2 } 
E_{\vartheta_2 \varphi_2 }(z_2 ,w_2 ) D_{z_2, \vartheta_2 , \varphi_2 } f_2 (z_2 )dx_2  dy_2 \edag,
\end{align*}
and the theorem follows.
\end{proof}
\begin{remark} Let us make the following comments:
\begin{enumerate}
\item If $\vartheta (Z)= \displaystyle\frac{1}{2} $ and $ \varphi(Z) =\displaystyle \frac{{\bf i}}{2} $  the previous results belong to the usual bicomplex holomorphic function theory. 
\item The function theory obtained  for  $\vartheta (Z)= \displaystyle\frac{1}{2}$ and $\varphi(Z) = \displaystyle -\frac{{\bf i}}{2} $  can be considered as the   bicomplex anti-holomorphic function theory. Set 
$$\frac{\partial}{\partial Z_{(\vartheta,\varphi)}} = \e\frac{\partial }{ \partial z_1}  + \edag \frac{\partial }{ \partial z_2}.$$
\item For $\vartheta (Z)= \displaystyle\frac{1}{2}$ and $\varphi(Z) =  \e  \displaystyle\frac{{\bf i}}{2} -   \edag \frac{{\bf i}}{2}$, set   
$$\frac{\partial}{\partial Z_{(\vartheta,\varphi)}} = \e\frac{\partial }{ \partial \bar z_1}  + \edag \frac{\partial }{ \partial z_2}.$$ Then we induce the space of functions $F=f_1\e +f_2 \edag$, where $f_1 $ is a holomorphic function in a complex variable and $f_2$ is an anti-holomorphic function. Similar type of behavior is presented for $\vartheta (Z)= \displaystyle\frac{1}{2}$ and $\varphi(Z) = \displaystyle-\frac{{\bf i}}{2}\e + \displaystyle\frac{{\bf i}}{2} \edag$.
\end{enumerate}
\end{remark}
\section{On a bicomplex $(\vartheta,\varphi)-$weighted type of\\ Riemann-Liouville fractional derivative}
Set ${U}=(a_1, c_1,a_2, c_2) , \ V=(b_1, d_1,b_2, d_2) \in \mathbb R^4$ such that  $a_l< b_l$  y $c_l < d_l$ for $l=1,2 $ and  denote 
\begin{align*}
J_U^V & = \{ z_1\e + z_2 \edag \in \bc \ \mid  \ \Re z_l \in [a_l,b_l], \ \Im z_l \in [c_l,d_l] , \  \ l=1,2  \} \\
 &= \left( [a_1,b_1]+ {\bf i} [c_1,d_1]\right) \e + \left( [a_2,b_2]+ {\bf i} [c_2,d_2]\right) \edag  
\end{align*}
In what follows, we write $\vec{\alpha}$ for vectors $(\alpha_0, \alpha_1,\alpha_2,\alpha_3)$ in $(0,1)^4$.
 
\begin{definition} \label{Def1}  Fixed  $W= w_1\e + w_2\edag \in J_{U}^V$ 
 and set   $F: J_{U}^V  \to \bc $. We shall write  $F\in AC^1(J_{U}^V , \bc)$ if there exist  
$f_l: [a_l,b_l]+ {\bf i} [c_l,d_l] \to \mathbb C({\bf i})$ for $l=1,2$ such that  the maps 
$ \Re z_l \to f_l (\Re z_l + {\bf i} \Im w_l)$ belong  to  $AC^1([a_l,b_l], \C({\bf i}))$ and  
   $ \Im z_l \to f_l (\Re w_l + {\bf i} \Im z_l)$ belong to $AC^1( [c_l,d_l], \C({\bf i}))$ for $l=1,2$. 
The right bicomplex $(\vartheta,\varphi)-$fractional derivative  in the Riemann-Liouville sense, or  bicomplex $(\vartheta,\varphi)-$fractional Riemann-Liouville derivative of $F$ is defined by 
\begin{align*} 
&( {}^{(\vartheta,\varphi)}\mathcal D_{a^+   }^{\vec \alpha} F)(W,Z):= \left( \vartheta_1 
(D _{a_1^+}^{\alpha_0} f_1)(x_1+ {\bf i}\Im w_1 )  +\varphi_1 
	(D _{c_1^+}^{\alpha_1} f_1)(\Re w_1 + {\bf i} y_1 )  \right) \e  \\
	 &  +  
	 \left( \vartheta_2 
(D _{a_2^+}^{\alpha_2} f_2)(x_2 + {\bf i} \Im w_2 )  +\varphi_2 
	(D _{c_2^+}^{\alpha_3} f_2)(\Re w_2 + {\bf i} y_2 )  \right) \edag \\
	= &   \left( \vartheta_1 
\frac{\partial}{\partial{x_1}} \left[ ({\bf I}_{a_0^+}^{1-\alpha_0} f_1)(x_1+ {\bf i}\Im w_1 ) \right]
  +\varphi_1 
	\frac{\partial}{\partial{y_1}} \left[ ({\bf I}_{a_1^+} ^{1-\alpha_1} f_1)(\Re w_1 + {\bf i} y_1 ) \right]
  \right) \e    \\
	& +  
	\left( \vartheta_2
		\frac{\partial}{\partial{x_2}} \left[ ({\bf I}_{a_2^+} ^{1-\alpha_2} f_2)(x_2 + {\bf i} \Im w_2 ) \right]
  +\varphi_2 
		\frac{\partial}{\partial{y_2}} \left[ ({\bf I}_{a_3^+} ^{1-\alpha_3} f_2)(\Re w_2 + {\bf i} y_2 )\right]
  \right) \edag   
\end{align*}
and its left version is given by 
 \begin{align*} 
( {}^{(\vartheta,\varphi)}\mathcal D_{b^-}^{\vec \alpha} F)(W,Z)  
:= &   \left( \vartheta_1 
(D _{b_1^-}^{\alpha_0} f_1)(x_1+ {\bf i}\Im w_1 )  +\varphi_1 
	(D _{d_1^-}^{\alpha_1} f_1)(\Re w_1 + {\bf i} y_1 )  \right) \e  \\
	 &  +  
	 \left( \vartheta_2 
(D _{b_2^-}^{\alpha_2} f_2)(x_2 + {\bf i} \Im w_2 )  +\varphi_2 
	(D _{d_2^-}^{\alpha_3} f_2)(\Re w_2 + {\bf i} y_2 )  \right) \edag.  
\end{align*}
Define the following operators:
\begin{align*}
 {}^{(\vartheta,\varphi)} \mathcal I_{a}^{\vec{\alpha}} [F](W,Z): =  & \left[   \left(  \  ({\bf I}_{a_0^+}^{1-\alpha_0} f_1)(x_1+ {\bf i}\Im w_1 ) 
  +  ({\bf I}_{a_1^+} ^{1-\alpha_1} f_1)(\Re w_1 + {\bf i} y_1 )   \right) \e \right. \\
		 &  \ \ \left.	+  
	\left(   \   ({\bf I}_{a_2^+} ^{1-\alpha_2} f_2)(x_2 + {\bf i} \Im w_2 ) + 
	({\bf I}_{a_3^+} ^{1-\alpha_3} f_2)(\Re w_2 + {\bf i} y_2 ) \right) \edag \right], 
\\
{}^{(\vartheta,\varphi)}\mathfrak P_a^{\vec \alpha} := &
   \left(  \  D_{a_0^+}^{1-\alpha_{0}}  
  +  D_{a_1^+}^{1-\alpha_{1}}    \right) \e   	+  
	\left(   \  D_{a_2^+}^{1-\alpha_{2}}    + 
	    D_{a_3^+}^{1-\alpha_{3}}\right) \edag,			
\end{align*}
where $D_{a_0^+}^{1-\alpha_{0}}, D_{a_1^+}^{1-\alpha_{1}}, D_{a_2^+}^{1-\alpha_{2}},  D_{a_3^+}^{1-\alpha_{3}}$ stand for the fractional partial derivatives with respect to the real variables $x_1, y_1, x_2, y_2$ respectively.
 
Considering $Z = (x_1+y_1{\bf i})\e + (x_2+y_2{\bf i})\edag$, we introduce 
\begin{align*}	
 & {}^{(\vartheta,\varphi)}\mathfrak T_a^{\vec \alpha} F(W,Z):=  \\
  &   \left(  \  \frac{(x_1-a_0)^{-\alpha_0}}{\Gamma[1-\alpha_0]}  ({\bf I}_{a_1^+} ^{1-\alpha_1} f_1)(\Re w_1 + {\bf i} y_1 )
  +  \frac{(y_1-a_1)^{-\alpha_1}}{\Gamma[1-\alpha_1]}  ({\bf I}_{a_0^+}^{1-\alpha_0} f_1)(x_1+ {\bf i}\Im w_1 )     \right) \e   \\
		 & \ \ 	+   
	\left(   \  \frac{(x_2-a_2)^{-\alpha_2}}{\Gamma[1-\alpha_2]}    ({\bf I}_{a_3^+} ^{1-\alpha_3} f_2)(\Re w_2 + {\bf i} y_2 ) 
	+  \frac{(y_2-a_3)^{-\alpha_3}}{\Gamma[1-\alpha_3]}   ({\bf I}_{a_2^+} ^{1-\alpha_2} f_2)(x_2 + {\bf i} \Im w_2 )   \right) \edag .   
\end{align*}
\end{definition}

\begin{proposition} \label{DI} According to Definition \ref{Def1} we have that  
\begin{align*} ({}^{(\vartheta,\varphi)}\mathcal D_{a^+   }^{\vec \alpha} F)(W,Z ) =   & \frac{\partial  }{\partial Z_{\vartheta \varphi } } \circ  {}^{(\vartheta,\varphi)} \mathcal I_{a}^{\vec{\alpha}} [F](W,Z), \\
   {}^{(\vartheta,\varphi)} \mathfrak P_a^{\vec \alpha} \circ   {}^{(\vartheta,\varphi)} \mathcal I_{a}^{\vec{\alpha}} [F](W,Z) 	= &   \left( f_1 (x_1+ {\bf i}\Im w_1 ) 
  +   f_1 (\Re w_1 + {\bf i} y_1 )   \right) \e \\ 
+ \left(f_2 (x_2 + {\bf i} \Im w_2 ) + f_2(\Re w_2 + {\bf i} y_2 ) \right)\edag  & +  {}^{(\vartheta,\varphi)} \mathfrak T_a^{\vec \alpha} F(W,Z).
\end{align*}
\end{proposition}
\begin{proof}
\begin{align*} 
&( {}^{(\vartheta,\varphi)}\mathcal D_{a^+   }^{\vec \alpha} F)(W,Z) \\
	= &   \left( \vartheta_1 
\frac{\partial }{\partial{x_1} } \left[ ({\bf I}_{a_0^+}^{1-\alpha_0} f_1)(x_1+ {\bf i}\Im w_1 ) \right]
  +\varphi_1 
	\frac{\partial }{\partial{y_1} } \left[ ({\bf I}_{a_1^+} ^{1-\alpha_1} f_1)(\Re w_1 + {\bf i} y_1 ) \right]
  \right) \e    \\
	& +  
	\left( \vartheta_2
		\frac{\partial }{\partial{x_2} } \left[ ({\bf I}_{a_2^+} ^{1-\alpha_2} f_2)(x_2 + {\bf i} \Im w_2 ) \right]
  +\varphi_2 
		\frac{\partial }{\partial{y_2} } \left[ ({\bf I}_{a_3^+} ^{1-\alpha_3} f_2)(\Re w_2 + {\bf i} y_2 )\right]
  \right) \edag   \\
	= &   \left(  \  ( \vartheta_1 
\frac{\partial }{\partial{x_1} }  
  +\varphi_1 
	\frac{\partial }{\partial{y_1} }     ) \ \e    +  
	(  \   \vartheta_2
		\frac{\partial }{\partial{x_2} }   +\varphi_2 
		\frac{\partial }{\partial{y_2} }    )  \  \edag \right)   \\
		&  \ \left[   \left(  \  ({\bf I}_{a_0^+}^{1-\alpha_0} f_1)(x_1+ {\bf i}\Im w_1 )  
  +  ({\bf I}_{a_1^+} ^{1-\alpha_1} f_1)(\Re w_1 + {\bf i} y_1 )   \right) \e   \right. 
		\\
		&  \ \  \  \  \left. +  
	\left(   \   ({\bf I}_{a_2^+} ^{1-\alpha_2} f_2)(x_2 + {\bf i} \Im w_2 ) + 
	({\bf I}_{a_3^+} ^{1-\alpha_3} f_2)(\Re w_2 + {\bf i} y_2 ) \right) \edag \right]
 \\
=& \frac{\partial  }{\partial Z_{\vartheta \varphi } } \circ  {}^{(\vartheta,\varphi)} \mathcal I_{a}^{\vec{\alpha}} [F](W,Z), \end{align*}
On the other hand, equations \eqref{FundTheorem} and \eqref{cte} imply that 
\begin{align*}
  & {}^{(\vartheta,\varphi)} \mathfrak P_a^{\vec \alpha} \circ   {}^{(\vartheta,\varphi)} \mathcal I_{a}^{\vec{\alpha}} [F](W,Z) \\
	= & 
		    \left(  \  D_{a_0^+}^{1-\alpha_{0}} ({\bf I}_{a_0^+}^{1-\alpha_0} f_1)(x_1+ {\bf i}\Im w_1 ) 
  +  D_{a_1^+}^{1-\alpha_{1}} ({\bf I}_{a_1^+} ^{1-\alpha_1} f_1)(\Re w_1 + {\bf i} y_1 )   \right) \e   \\
		 &  	+  
	\left(   \  D_{a_2^+}^{1-\alpha_{2}}  ({\bf I}_{a_2^+} ^{1-\alpha_2} f_2)(x_2 + {\bf i} \Im w_2 ) + 
	    D_{a_3^+}^{1-\alpha_{3}}  ({\bf I}_{a_3^+} ^{1-\alpha_3} f_2)(\Re w_2 + {\bf i} y_2 ) \right) \edag  \\ 
	   & +  
		    \left(  \  D_{a_0^+}^{1-\alpha_{0}}  ({\bf I}_{a_1^+} ^{1-\alpha_1} f_1)(\Re w_1 + {\bf i} y_1 )
  +  D_{a_1^+}^{1-\alpha_{1}}  ({\bf I}_{a_0^+}^{1-\alpha_0} f_1)(x_1+ {\bf i}\Im w_1 )     \right) \e   \\
		 &  	+  
	\left(   \  D_{a_2^+}^{1-\alpha_{2}}  ({\bf I}_{a_3^+} ^{1-\alpha_3} f_2)(\Re w_2 + {\bf i} y_2 ) + 
	    D_{a_3^+}^{1-\alpha_{3}}  ({\bf I}_{a_2^+} ^{1-\alpha_2} f_2)(x_2 + {\bf i} \Im w_2 )   \right) \edag  \\ 
=		&   \left( f_1 (x_1+ {\bf i}\Im w_1 ) 
  +   f_1 (\Re w_1 + {\bf i} y_1 )   \right) \e   + 
	\left(    f_2 (x_2 + {\bf i} \Im w_2 ) + 
	   f_2(\Re w_2 + {\bf i} y_2 ) \right) \edag \\
	   & +  
		    \left(  \  \frac{(x_1-a_0)^{-\alpha_0}}{\Gamma[1-\alpha_0]}  ({\bf I}_{a_1^+} ^{1-\alpha_1} f_1)(\Re w_1 + {\bf i} y_1 )
  +  \frac{(y_1-a_1)^{-\alpha_1}}{\Gamma[1-\alpha_1]}  ({\bf I}_{a_0^+}^{1-\alpha_0} f_1)(x_1+ {\bf i}\Im w_1 )     \right) \e   \\
		 &  	+  
	\left(   \  \frac{(x_2-a_2)^{-\alpha_2}}{\Gamma[1-\alpha_2]}    ({\bf I}_{a_3^+} ^{1-\alpha_3} f_2)(\Re w_2 + {\bf i} y_2 ) 
	+  \frac{(y_2-a_3)^{-\alpha_3}}{\Gamma[1-\alpha_3]}   ({\bf I}_{a_2^+} ^{1-\alpha_2} f_2)(x_2 + {\bf i} \Im w_2 )   \right) \edag.    
\end{align*}
 \end{proof}

\section{Main Results}
In this section we state and prove the main theoretical paper's results on bicomplex functions associated to ${}^{(\vartheta,\varphi)}\mathcal D_{a^+}^{\vec \alpha}$. Specifically, something like Gauss's theorem, Cauchy's theorem and a Borel-Pompieu formula as well.  

\begin{theorem} \label{GT}
(Gauss's theorem for bicomplex functions associated to ${}^{(\vartheta,\varphi)}\mathcal D_{a^+   }^{\vec \alpha} $)
Set $W=w_1\e+ w_2\edag \in \Gamma$ fixed, denote our bicomplex variable by $Z=z_1\e + z_2\edag \in J_U^V$ and set a bicomplex function $F\in AC^1(J_{U}^V , \bc)$. Then  
\begin{align*}  &  \int_{J_U^V} \left[({}^{(\vartheta,\varphi)}\mathcal D_{a^+}^{\vec \alpha} F)(W,Z)
 +  A_{\vartheta  \varphi }  {}^{(\vartheta,\varphi)} \mathcal I_{a}^{\vec{\alpha}} [F](W,Z)  + B_{\vartheta  \varphi }  {\bf i} {}^{(\vartheta,\varphi)} \mathcal I_{a}^{\vec{\alpha}} [F](W,Z)  \right] dZ \wedge dZ^{*} \\
 = &
\int_{\partial J_U^V}
{}^{(\vartheta,\varphi)} \mathcal I_{a}^{\vec{\alpha}} [F](W,Z) d\sigma_{(\vartheta , \varphi) }(Z),
\end{align*}
\end{theorem}
\begin{proof}
Applying Theorem \ref{weighted_Gauss_Bicomplex} for ${}^{(\vartheta,\varphi)} \mathcal I_{a}^{\vec{\alpha}} [F](W,Z)$ and by Proposition \ref{DI} the result follows.
\end{proof}

\begin{corollary} 
(Cauchy's theorem for bicomplex functions associated to ${}^{(\vartheta,\varphi)}\mathcal D_{a^+}^{\vec \alpha} $)  According to Theorem \ref{GT}  and if $( {}^{(\vartheta,\varphi)}\mathcal D_{a^+}^{\vec \alpha} F)(W,Z)=0$ for all $Z\in \Gamma$ then  
\begin{align*}  &  \int_{J_U^V} \left(A_{\vartheta  \varphi }  {}^{(\vartheta,\varphi)} \mathcal I_{a}^{\vec{\alpha}} [F](W,Z)  + B_{\vartheta  \varphi }  {\bf i} {}^{(\vartheta,\varphi)} \mathcal I_{a}^{\vec{\alpha}} [F](W,Z)  \right) dZ \wedge dZ^{*} \\
 = &
\int_{\partial J_U^V}
{}^{(\vartheta,\varphi)} \mathcal I_{a}^{\vec{\alpha}} [F](W,Z) d\sigma_{(\vartheta ,\varphi )}(Z) ,
\end{align*}
\end{corollary}
We are in conditions to establish a weighted version of the Borel-Pompeiu formula in a bicomplex fractional setting.

\begin{theorem}\label{BPF} (Borel-Pompieu formula associated to ${}^{(\vartheta,\varphi)}\mathcal D_{a^+}^{\vec \alpha}$) 
Set $W=w_1\e+ w_2\edag  , Q=(q_1+{\bf i} r_1)\e + (q_2+{\bf i} r_2)\edag  \in \Gamma$, denote our bicomplex variable by $Z=z_1\e + z_2\edag \in J_U^V$ and set  $F=f_1\e +f_2\edag \in AC^1(J_{U}^V , \bc)$. Given  $A\in J_{U}^V$  and $r\in \mathbb D^+$ such that 
 $\overline{B(A, r)} \subset J_{U}^V $ and fixed $W =w_1\e + w_2 \edag  \in J_{U}^V$, where $w_1,w_2$ are complex numbers. 
Then, we have that 
\begin{align*}  
&   \left( f_1 (q_1+ {\bf i}\Im w_1 ) 
  +   f_1 (\Re w_1 + {\bf i} r_1 )   \right) \e    + 
	\left(    f_2 (q_2 + {\bf i} \Im w_2 ) + 
	   f_2(\Re w_2 + {\bf i} r_2 ) \right) \edag\\
	  =   & 
\int_{\partial B(A, r) }  
{}^{(\vartheta,\varphi)} \mathcal I_{a}^{\vec{\alpha}} [F](W,Z)  \ {}^{(\vartheta,\varphi)}\mathfrak P_a^{\vec \alpha} E_{(\vartheta, \varphi)} (Z , Q )d\sigma_{(\vartheta, \varphi)}(Z ) 
 \\
 &  \   \  - \int_{B(A, r)} \ {}^{(\vartheta,\varphi)}\mathfrak P_a^{\vec \alpha} E_{(\vartheta, \varphi)}(Z ,Q ) ( {}^{(\vartheta,\varphi)}\mathcal D_{a^+   }^{\vec \alpha} F)(W,Z) dZ\wedge dZ^{*}   -  {}^{(\vartheta,\varphi)} \mathfrak T_a^{\vec \alpha} F(W,Q).
\end{align*}
for every $Q = (q_1+ r_1)\e + (q_2+ r_2)\edag  \in B(A, r)$, where $q_1,r_1,q_2,r_2$ are real numbers.
\end{theorem}

\begin{proof} 
Applying Theorem \ref{weighted_Borel-Pompieu_Bicomplex} for ${}^{(\vartheta,\varphi)} \mathcal I_{a}^{\vec{\alpha}} [F](W,Z)$ and by Proposition \ref{DI} we get that  
\begin{align*}   
{}^{(\vartheta,\varphi)} \mathcal I_{a}^{\vec{\alpha}} [F](W,Q) c_{(\vartheta,\varphi)}   & =
\int_{\partial  B(A, r)  }  
{}^{(\vartheta,\varphi)} \mathcal I_{a}^{\vec{\alpha}} [F](W,Z)  E_{(\vartheta, \varphi)} (Z , Q )d\sigma_{(\vartheta, \varphi)}(Z ) 
 \\
 &  \   \  - \int_{ B(A, r)  } 
E_{\vartheta , \varphi }(Z ,Q ) ( {}^{(\vartheta,\varphi)}\mathcal D_{a^+}^{\vec \alpha} F)(W,Z) dZ\wedge dZ^{*}.
\end{align*}

We now apply operator ${}^{(\vartheta,\varphi)}\mathfrak P_a^{\vec \alpha} $ to both sides, with respect to  the real variables of $Q$ (the second component of domain's elements of $E_{(\vartheta, \varphi)}(Z,Q)$).

From the  properties of the complex  weighted Cauchy kernel  $$ E_{(\vartheta_1, \varphi_1)} (z_1 , q_1 )= \sum_{n=0}^{\infty} A_n(z_1)C_n(q_1 ), \ 
E_{(\vartheta_ 2 , \varphi_ 2 )} (z_ 2  , q_ 2  )= \sum_{n=0}^{\infty} B_n(z_ 2 ) C_n(q_ 2  ),$$
we have that 
$$E_{(\vartheta, \varphi)} (Z,Q) = \sum_{n=0}^{\infty} \left( A_n(z_1)  {\bf e }  + B_n(z_ 2 ) {\bf e }^{\dagger} \right)      \left( C_n(q_1)  {\bf e }  + D_n(q_ 2 ) {\bf e }^{\dagger} \right)= \sum_{n=0}^{\infty} U_n(Z) V_n(Q),$$
where the convergence is  uniform on compact subsets of $\mathbb {BC}$ with  $Z\neq Q$.  
 
Therefore, using the uniform convergence of the previous series we have that  
\begin{align*} & {}^{(\vartheta,\varphi)}\mathfrak P_a^{\vec \alpha} \int_{\partial   B(A, r)  }  
{}^{(\vartheta,\varphi)} \mathcal I_{a}^{\vec{\alpha}} [F](W,Z)  E_{(\vartheta, \varphi)} (Z , Q )d\sigma_{(\vartheta, \varphi)}(Z ) \\    
= &  {}^{(\vartheta,\varphi)}\mathfrak P_a^{\vec \alpha} \int_{\partial  B(A, r)  }  
{}^{(\vartheta,\varphi)} \mathcal I_{a}^{\vec{\alpha}} [F](W,Z)  \left(\sum_{n=0}^{\infty} U_n(Z) V_n(Q) \right)d\sigma_{(\vartheta, \varphi)}(Z ) \\
= & \sum_{n=0}^{\infty}   \int_{  \partial B(A, r)  }  
{}^{(\vartheta,\varphi)} \mathcal I_{a}^{\vec{\alpha}} [F](W,Z)   U_n(Z) {}^{(\vartheta,\varphi)}\mathfrak P_a^{\vec \alpha} (V_n(Q) )   d\sigma_{(\vartheta, \varphi)}(Z )   \\
= &   \int_{\partial   B(A, r)  }  
{}^{(\vartheta,\varphi)} \mathcal I_{a}^{\vec{\alpha}} [F](W,Z) \sum_{n=0}^{\infty} \left( U_n(Z) {}^{(\vartheta,\varphi)}\mathfrak P_a^{\vec \alpha} (V_n(Q) ) \right)  d\sigma_{(\vartheta, \varphi)}(Z ) \\
 =&  \int_{\partial  B(A, r)  }  
{}^{(\vartheta,\varphi)} \mathcal I_{a}^{\vec{\alpha}} [F](W,Z) {}^{(\vartheta,\varphi)}\mathfrak P_a^{\vec \alpha} E_{(\vartheta, \varphi)} (Z , Q )d\sigma_{(\vartheta, \varphi)}(Z).
\end{align*}
Note that the equality
 \begin{align*} &  {}^{(\vartheta,\varphi)} \mathfrak P_a^{\vec \alpha} \int_{  B(A, r)} 
   E_{(\vartheta, \varphi)}(Z ,Q ) ( {}^{(\vartheta,\varphi)}\mathcal D_{a^+}^{\vec \alpha} F)(W,Z) dZ\wedge dZ^{*} \\
	=  & \int_{  B(A, r)} 
 \ {}^{(\vartheta,\varphi)} \mathfrak P_a^{\vec \alpha} E_{(\vartheta, \varphi)}(Z ,Q ) ( {}^{(\vartheta,\varphi)}\mathcal D_{a^+}^{\vec \alpha} F)(W,Z) dZ\wedge dZ^{*}
\end{align*}
is achieved by Fubini's Theorem (the switching of the order of the iterated integrals) and Leibniz's rule. Consequently,  
\begin{align*}   
& {}^{(\vartheta,\varphi)} \mathfrak P_a^{\vec \alpha} \circ {}^{(\vartheta,\varphi)} \mathcal I_{a}^{\vec{\alpha}} [F]( W,Q) c_{(\vartheta, \varphi)}   \\
 = &
\int_{\partial  B(A, r)}  
{}^{(\vartheta,\varphi)} \mathcal I_{a}^{\vec{\alpha}} [F](W,Z)  \  {}^{(\vartheta,\varphi)}\mathfrak P_a^{\vec \alpha} E_{(\vartheta, \varphi)} (Z , Q )d\sigma_{(\vartheta, \varphi)}(Z ) 
 \\
 &  \   \  - \int_{  B(A, r)} 
 \ {}^{(\vartheta,\varphi)} \mathfrak P_a^{\vec \alpha} E_{(\vartheta, \varphi)}(Z ,Q ) ( {}^{(\vartheta,\varphi)}\mathcal D_{a^+}^{\vec \alpha} F)(W,Z) dZ\wedge dZ^{*}
\end{align*}
and the proof is complete from Proposition \ref{DI}.
\end{proof}

\begin{corollary}
(Cauchy formula associated to ${}^{(\vartheta,\varphi)}\mathcal D_{a^+}^{\vec \alpha} $) 
 Under the same hypothesis of Theorem \ref{BPF}, if moreover $({}^{(\vartheta,\varphi)}\mathcal D_{a^+}^{\vec \alpha} F)(W,Z)=0$ for every $Z\in \Gamma$, then 
\begin{align*}  
&   \left( f_1 (q_1+ {\bf i}\Im w_1 ) 
  +   f_1 (\Re w_1 + {\bf i} r_1 )   \right) \e    + 
	\left(    f_2 (q_2 + {\bf i} \Im w_2 ) + 
	   f_2(\Re w_2 + {\bf i} r_2 ) \right) \edag\\
	  =   & 
\int_{\partial   B(A, r) }  
{}^{(\vartheta,\varphi)} \mathcal I_{a}^{\vec{\alpha}} [F](W,Z)   \  {}^{(\vartheta,\varphi)} \mathfrak P_a^{\vec \alpha} E_{(\vartheta, \varphi)} (Z , Q)d\sigma_{(\vartheta, \varphi)}(Z )  - {}^{(\vartheta,\varphi)} \mathfrak T_a^{\vec \alpha} F(W,Q).
\end{align*}\end{corollary}

\begin{remark} Here are some particular situations in which we can draw conclusions for much traditional theories
\begin{enumerate}
\item If  $\vartheta (Z)= 1  $ and $ \varphi(Z) = {\bf i}$  the previous results describe the behavior of a usual bicomplex holomorphic fractional derivative. 
\item For  $\vartheta (Z)= 1  $ and $ \varphi(Z) = -{\bf i}$  we obtain the properties of a bicomplex anti-holomorphic fractional derivative. 
\item Suppose that  $\vartheta (Z)= 1   $ and $ \varphi(Z) = {\bf i}\e - {\bf i} \edag$ then the component of $\e$ is reduced to  holomorphic fractional derivative, meanwhile in $\edag$  appeared anti-holomorphic fractional derivative. Similar type of behavior is presented for $\vartheta (Z)= 1   $ and $ \varphi(Z) = -{\bf i}\e + {\bf i} \edag$.
\end{enumerate}
\end{remark}
 
\section*{Acknowledgments}
This work was partially supported by the Instituto Polit\'ecnico Nacional (grant numbers SIP20220017, SIP20221274).

\end{document}